\documentclass{amsart}
\usepackage{amsmath}
\usepackage[dvips]{graphicx}
\usepackage{verbatim}
\date{\today}
\allowdisplaybreaks
\theoremstyle{plain}
  \newtheorem{thm}{Theorem}
  \newtheorem{lem}{Lemma}[section]

\theoremstyle{definition}

\theoremstyle{remark}
  
  \newtheorem*{ack}{Acknowledgment}
\numberwithin{equation}{section}
\newcommand{\Z}{{\mathbf{Z}}}
\newcommand{\bk}{{\mathbf{k}}}
\newcommand{\cL}{{\mathcal{L}}}
\newcommand{\cO}{{\mathcal{O}}}
\newcommand{\lang}{\left\langle}
\newcommand{\rang}{\right\rangle}
\thicklines
\newcommand{\posi}{
  \raisebox{-4mm}{
  \begin{picture}(11.5,10)
  \put( 0, 0){\line( 1, 1){10}}
  \put( 0,10){\line( 1,-1){ 4}}
  \put(10, 0){\line(-1, 1){ 4}}
  \end{picture}}}
\newcommand{\zero}{
  \raisebox{-4mm}{
  \begin{picture}(11.5,10)
  \qbezier( 0, 0)( 4, 5)( 0,10)
  \qbezier(10, 0)( 6, 5)(10,10)
  \end{picture}}}
\newcommand{\infi}{
  \raisebox{-4mm}{
  \begin{picture}(11.5,10)
  \qbezier( 0, 0)( 5, 4)(10, 0)
  \qbezier( 0,10)( 5, 6)(10,10)
  \end{picture}}}
\newcommand{\posneg}{
  \raisebox{-4mm}{
  \begin{picture}(26.5,16)
  \put( 0, 1){\line( 1, 1){ 9}}
  \put( 1, 0){\line( 1, 1){ 9}}
  \put( 0,10){\line( 1,-1){ 4}}
  \put(10, 0){\line(-1, 1){ 4}}
  \put(16,10){\line( 1,-1){ 9}}
  \put(15, 9){\line( 1,-1){ 9}}
  \put(15, 0){\line( 1, 1){ 4}}
  \put(25,10){\line(-1,-1){ 4}}
  \thinlines
  \qbezier[10](10, 9)(12.5,12.5)(15, 9)
  \qbezier[10]( 9,10)(12.5,13.5)(16,10)
  \qbezier[10](10, 0)(12.5,-2.5)(15, 0)
  \qbezier[25]( 0,10)(-2,12)(12.5,14)
  \qbezier[25](25,10)(27,12)(12.5,14)
  \qbezier[25]( 1, 0)(-1,-2)(12.5,-3)
  \qbezier[25](24, 0)(26,-2)(12.5,-3)
  \qbezier[25]( 0, 1)(-4,-3)(12.5,-4)
  \qbezier[25](25, 1)(29,-3)(12.5,-4)
  \end{picture}}}
\newcommand{\negpos}{
  \raisebox{-4mm}{
  \begin{picture}(26.5,16)
  \put( 0, 1){\line( 1, 1){ 4}}
  \put( 1, 0){\line( 1, 1){ 4}}
  \put( 9,10){\line(-1,-1){ 4}}
  \put(10, 9){\line(-1,-1){ 4}}
  \put( 0,10){\line( 1,-1){10}}
  \put(16,10){\line( 1,-1){ 4}}
  \put(15, 9){\line( 1,-1){ 4}}
  \put(25, 1){\line(-1, 1){ 4}}
  \put(24, 0){\line(-1, 1){ 4}}
  \put(15, 0){\line( 1, 1){10}}
  \thinlines
  \qbezier[10](10, 9)(12.5,12.5)(15, 9)
  \qbezier[10]( 9,10)(12.5,13.5)(16,10)
  \qbezier[10](10, 0)(12.5,-2.5)(15, 0)
  \qbezier[25]( 0,10)(-2,12)(12.5,14)
  \qbezier[25](25,10)(27,12)(12.5,14)
  \qbezier[25]( 1, 0)(-1,-2)(12.5,-3)
  \qbezier[25](24, 0)(26,-2)(12.5,-3)
  \qbezier[25]( 0, 1)(-4,-3)(12.5,-4)
  \qbezier[25](25, 1)(29,-3)(12.5,-4)
  \end{picture}}}
\newcommand{\nulinf}{
  \raisebox{-4mm}{
  \begin{picture}(26.5,16)
  \qbezier( 0, 1)( 5, 5)( 0,10)
  \qbezier(10, 0)( 5, 5)(10, 9)
  \qbezier( 1, 0)( 5, 3)( 5, 5)
  \qbezier( 9,10)( 5, 7)( 5, 5)
  \qbezier(16,10)(20, 6)(25,10)
  \qbezier(15, 0)(20, 4)(24, 0)
  \qbezier(15, 9)(18, 5)(20, 5)
  \qbezier(25, 1)(22, 5)(20, 5)
  \thinlines
  \qbezier[10](10, 9)(12.5,12.5)(15, 9)
  \qbezier[10]( 9,10)(12.5,13.5)(16,10)
  \qbezier[10](10, 0)(12.5,-2.5)(15, 0)
  \qbezier[25]( 0,10)(-2,12)(12.5,14)
  \qbezier[25](25,10)(27,12)(12.5,14)
  \qbezier[25]( 1, 0)(-1,-2)(12.5,-3)
  \qbezier[25](24, 0)(26,-2)(12.5,-3)
  \qbezier[25]( 0, 1)(-4,-3)(12.5,-4)
  \qbezier[25](25, 1)(29,-3)(12.5,-4)
  \end{picture}}}
\newcommand{\infnul}{
  \raisebox{-4mm}{
  \begin{picture}(26.5,16)
  \qbezier( 1, 0)( 5, 5)(10, 0)
  \qbezier( 0,10)( 5, 5)( 9,10)
  \qbezier( 0, 1)( 3, 5)( 5, 5)
  \qbezier(10, 9)( 7, 5)( 5, 5)
  \qbezier(15, 9)(20, 5)(15, 0)
  \qbezier(25, 1)(20, 5)(25,10)
  \qbezier(16,10)(20, 7)(20, 5)
  \qbezier(24, 0)(20, 3)(20, 5)
  \thinlines
  \qbezier[10](10, 9)(12.5,12.5)(15, 9)
  \qbezier[10]( 9,10)(12.5,13.5)(16,10)
  \qbezier[10](10, 0)(12.5,-2.5)(15, 0)
  \qbezier[25]( 0,10)(-2,12)(12.5,14)
  \qbezier[25](25,10)(27,12)(12.5,14)
  \qbezier[25]( 1, 0)(-1,-2)(12.5,-3)
  \qbezier[25](24, 0)(26,-2)(12.5,-3)
  \qbezier[25]( 0, 1)(-4,-3)(12.5,-4)
  \qbezier[25](25, 1)(29,-3)(12.5,-4)
  \end{picture}}}

\begin{document}
\title[$SU(2)$-invariants associated with non-trivial cohomology
classes]
{Quantum $SU(2)$-invariants for three-manifolds associated with
non-trivial cohomology classes modulo two}
\author{Hitoshi Murakami}
\address{Department of Mathematics, School of Science and Engineering,
Waseda University, Ohkubo, Shinjuku, TOKYO 169-8555, JAPAN,
and
Mittag-Leffler Institute, Aurav{\"a}gen 17, S-182 62,
Djursholm, SWEDEN}
\email{hitoshi@uguisu.co.jp}
\begin{abstract}
We show an integrality of the quantum $SU(2)$-invariant associated with a non-trivial
first cohomology class modulo two.
\end{abstract}
\thanks{
Partially supported by Waseda University Grant
for Special Research Projects (No. 98A-623) and
Grant-in-Aid for Scientific
Research (C) (No. 09640135), the Ministry of Education, Science, Sports and
Culture.}
\maketitle
For a given Lie group $G$ and an integer (level) $k$, E.~Witten gave an `invariant' of
three-manifolds by using the path integral \cite{Witten:COMMP89}.
Mathematically rigorous proof of its existence was given by several people.
See for example
\cite{Reshetikhin/Turaev:INVEM91,
      Kirby/Melvin:INVEM91,
      Blanchet/Habegger/Masbaum/Vogel:invariant_from_Kauffman_bracketTOPOL92,
      Wenzl:INVEM93,
      Turaev/Wenzl:INTJM693}.
There are also some refinements for these invariants corresponding to various structures
of the three-manifolds, such as cohomology classes and spin structures
\cite{Turaev:ICM90,
      Kirby/Melvin:INVEM91,
      Kohno/Takata:JKNOT93,
      Murakami:quantum,
      Blanchet:Warsaw,
      Blanchet:1998}.
\par
Some topological properties of the invariants were given especially for $G=SU(2)$.
The first striking result was found by R.~Kirby and P.~Melvin stating that the 
quantum $SU(2)$-invariant with level two splits into the sum of the invariants associated
with spin structures and each summand can be described in terms of the $\mu$-invariant.
They and X.~Zhang obtained a topological interpretation for the quantum $SU(2)$-invariant
with level four \cite{Kirby/Melvin/Zhang:COMMP93}.
In this case the invariant splits into the sum of the invariants associated with the
first cohomology classes modulo two and each summand can be described in terms of the
Witt invariant and the cohomologies of both the three-manifold and its double cover
defined by the cohomology class.
Similarly the author gave a topological interpretation for the quantum $SU(3)$-invariant
with level three \cite{Murakami:quantum}.
It also decomposes into the sum of the invariants associated with the first cohomology
classes modulo three and each summand can be described in terms of the cohomologies of
the manifold and its triple cover.
\par
If $G=SU(n)$ and $n$ and the level $k$ are coprime, then the invariant splits into the
product of the $U(1)$-invariant \cite{Murakami/Ohtsuki/Okada:ORSJM92,Gocho:PROJAA90}
and $PSU(n)$-invariant \cite{Kirby/Melvin:INVEM91,Kohno/Takata:1996}.
If $n+k$ is odd prime and the three-manifold is a rational homology sphere, the invariant
is a cyclotomic integer
\cite{Murakami:MATPC93,
      Murakami:MATPC95,
      Masbaum/Roberts:MATPC97,
      Masbaum/Wenzl:1997,
      Takata/Yokota:1996}.
Similar integrality holds for a three-manifold with positive first Betti number.
Moreover one can define the perturbative invariant from the $PSU(n)$-invariant for
rational homology three-spheres
\cite{Ohtsuki:MATPC95,
      Ohtsuki:INVEM96,
      Le:1998}.
Recently it is proved by T.T.Q.~Le that every coefficient of the perturbative invariant
for integral homology three-spheres is of finite type in T.~Ohtsuki's sense
\cite{Ohtsuki:JKNOT96}.
(For $n=2$ case, it is proved in by A.~Kricker and B.~Spence \cite{Kricker/Spence:JKNOT97}.
In particular its first non-trivial coefficient coincides with (a multiple of)
the Casson invariant \cite{Akbulut/McCarthy:Casson_invariant}.
(For rational homology three-spheres, Le showed in \cite{Le:1998} that this
coincides with a multiple of the Casson-Walker invariant
\cite{Walker:Casson_invariant}.)
See also \cite{Ohtsuki:1998} for recent developments of these invariants; especially
their relation with the LMO invariant \cite{Le/Murakami/Ohtsuki:TOPOL98}.
\par
In this paper we consider the quantum $SU(2)$-invariant with level divisible by four.
By  V.~Turaev, and Kirby and Melvin \cite{Kirby/Melvin:INVEM91,Turaev:ICM90} it
decomposes into the sum of the invariants associated with the cohomology classes modulo
two.
We will consider the case where $(n+2)/2$ is an odd prime.
\par 
Let $p$ be an odd prime and $M$ a $\Z/p\Z$-homology three-sphere.
For a cohomology class $\theta\in H^1(M;\Z/2\Z)$,
let $\tau_{2p}^{SU(2)}(M,\theta)$ be the quantum $SU(2)$-invariant with level $2p-2$
associated with the $2p$th root of unity $q=\exp(\pi\sqrt{-1}/p)$.
In the previous paper \cite{Murakami:1998} the author proved
\begin{thm}
  $\tau_{2p}^{SU(2)}(M,0)\in\Z[1/2,\xi]$
and the `coefficient' of $(\xi-1)$ is congruent to the Casson-Walker invariant
modulo $p$, where $\xi=\exp(2\pi\sqrt{-1}/p)$.
\end{thm}
\par
In this paper we will show a weaker result for non-trivial cohomology classes.
We will show
\begin{thm}\label{thm:main}
  \begin{equation*}
  \tau_{2p}^{SU(2)}(M,\theta)\in
  \begin{cases}
    \phantom{\sqrt{-1}}(\xi-1)\Z[1/2,\xi]\qquad\text{if $\theta\cup\theta\cup\theta = 0$},
    \\
             \sqrt{-1} (\xi-1)\Z[1/2,\xi]\qquad\text{if $\theta\cup\theta\cup\theta\ne0$}.
  \end{cases}
  \end{equation*}
\end{thm}
\begin{ack}
The author thanks the organizers, especially Sofia Lambropoulou, of the conference
`Knots in Hellas' for their hospitality.
\end{ack}
\section{proof}
We put $p=4u+\varepsilon$ with $\varepsilon=\pm1$.
Put $\xi=\exp(2\pi\sqrt{-1}/p)$, $q=\exp(\pi\sqrt{-1}/p)$ and
$s=\exp(\pi\sqrt{-1}/2p)$.
Therefore $\xi=q^2=s^4$.
\par
Suppose that a closed three-manifold $M$ is given by an algebraically split,
framed link $\cL$ with $m$ components.
We also assume that the framing of $\cL$ is $(f_1,f_2,\dots,f_m)$
with $f_i\not\equiv0\mod p$.
Moreover we assume that $f_i\equiv0\mod 2$ for $i\le b$ and $f_i\equiv1\mod 2$ for $b<i\le m$.
Let $\theta\in H^1(M;\Z/2\Z)$ be presented by the sublink
$\cL_{\theta}=L_1\cup L_2\cup\dots\cup L_c$ ($c\le b$).
We put $f_i=2g_i$ for $1\le i\le b$.
\par
By \cite[(8.32) Theorem]{Kirby/Melvin:INVEM91}, $\tau^{SU(2)}_{2p}(M,\theta)
=\alpha(\cL)\Sigma(\cL)$ where $\alpha(\cL)$ and $\Sigma(\cL)$ is given as
follows.
\begin{align*}
  \alpha(\cL)&=\left(\sqrt{\frac{1}{p}}\sin\left(\frac{\pi}{2p}\right)\right)^{m}
              \exp\left(-\frac{3(p-1)}{8p}\times2\pi\sqrt{-1}\right)^{\sigma(\cL)}
  \\
  \Sigma(\cL)&=\sum_{\substack{0<k_i\le p \\ \cL_{\bk}=\cL_{\theta}}}
              2^{|T_{\bk}|}[\bk]J(\cL,\bk),
\end{align*}
where $J(\cL,\bk)$ is the colored Jones polynomial of $\cL$ with color
$\bk=(k_1,k_2,\dots,k_m)$ ($0<k_i\le p$), $\cL_{\bk}$ is the sublink of $\cL$ consisting of
components with the color $k_i$ even, $|T_{\bk}|$ is the number of $k_i$ which
are less than $p$ and $\sigma(\cL)$ is the number of positives minus that of negatives
in $\{f_1,f_2,\dots f_m\}$.
\par
After some calculations we can simplify $\alpha(\cL)$ and $\Sigma(\cL)$ as
follows.
\begin{lem}[Lemma 2.2 in \cite{Murakami:1998}]
  \begin{equation*}
    \alpha(\cL)
    =
    \frac{(\xi^{u}+\xi^{-u})^m}{(2G(\xi))^m}
    (-1)^{u\sigma(\cL)+\sigma_{-}(\cL)(\varepsilon-1)/2}
    \xi^{3\times\overline{8}\varepsilon\sigma(\cL)},
  \end{equation*}
where $G(\xi)$ is a Gaussian sum $\sum_{k=0}^{p-1}\xi^{k^2}$, $\overline{8}$
is the inverse of $8$ in $\Z/p\Z$ and $\sigma_{-}(\cL)$ is the number of negatives
in  $\{f_1,f_2,\dots f_m\}$.
Note that $G(\xi)$ is of the form $\gamma(\xi-1)^{(p-1)/2}$ with $\gamma$ invertible
in $\Z[\xi]$.
\end{lem}
\begin{lem}[see Lemma 2.3 in \cite{Murakami:1998}]
  \begin{align*}
    \Sigma(\cL)&=[2]2^c\sum_{n_1,n_2,\dots,n_m=1}^{(p-1)/2}(-1)^{\sum_i n_i}
    V(L^{(2n_1-1,2n_2-1,\dots,2n_c-1,2n_{c+1},2n_{c+2},\dots,2n_{m})})
    \\
    &\times
    \prod_{l=1}^{c}C_{n_l}(g_l)
    \prod_{l=c+1}^{m}S_{n_l}(f_l)
  \\
  \intertext{with}
    S_{n}(f)&=\sum_{k=n}^{p-n-1}(-1)^{k}[2k+1]{\xi}^{(k^2+k)f/2}\binom{k+n}{k-n}
  \\
  \intertext{and}
    C_{n}(g)&=\sum_{k=n}^{(p-1)/2}(-1)^{k}[2k]s^{(4k^2-1)g}\binom{k+n-1}{k-n},
  \end{align*}
where $L^{(d_1,d_2,\dots,d_m)}$ is the (unframed) link obtained from $\cL$ by
replacing $i$th component with $d_i$ parallels respecting the zero-framing.
For a link $L$, $V(L)$ is the Jones polynomial defined by the skein relation
$q V(L_{+})-q^{-1}V(L_{-})=(s-s^{-1})V(L_{0})$ for the usual skein triple
$(L_{+},L_{-},L_{0})$ and by $V(\text{trivial knot})=1$.
\end{lem}
\begin{proof}
The proof is similar to \cite[Lemma 2.3]{Murakami:1998} and so details are omitted.
\end{proof}
\par
To prove Theorem~\ref{thm:main}, we prepare three lemmas.
\begin{lem}[Lemma 2.4 in \cite{Murakami:1998}]
  $S_{n}(f)\in\Z[\xi]$ is divisible by $(\xi-1)^{(p-1)/2-n}$ in $\Z[\xi]$.
\end{lem}
\begin{lem}
  $\sqrt{-1}^{g+1}C_{n}(g)\in\Z[\xi]$ is divisible by $(\xi-1)^{(p-1)/2-n+1}$
  in $\Z[\xi]$.
\end{lem}
\begin{proof}
First note that $s=\sqrt{-1}^{\varepsilon}\xi^{-\varepsilon u}$.
Therefore an easily calculation shows that
\begin{align*}
  C_n(g)=
  \frac{\sqrt{-1}^{-\varepsilon(g+1)}\xi^{\varepsilon ug}}
       {\xi^{u}+\xi^{-u}}
  \sum_{k=n}^{(p-1)/2}(\xi^{-2\varepsilon uk}-\xi^{2\varepsilon uk})
     \xi^{gk^2}\binom{k+n-1}{k-n}.
\end{align*}
Now the proof follows from \cite[Lemma 6.2]{Lin/Wang:95} (see also
\cite{Ohtsuki:MATPC95}).
\end{proof}
The following lemma is a generalization of a result of M.~Sakuma
\cite[Theorems 1 and 2]{Sakuma:MATPC88}, who proves the case where
$c=0$.
\begin{lem}
  $\sqrt{-1}^{c+1}
  V(L^{(2n_1-1,2n_2-1,\dots,2n_c-1,2n_{c+1},2n_{c+2},\dots,2n_{m})})
  \in\Z[\xi]$
  is divisible by $(\xi-1)^{\sum_{l=1}^{m}n_l-c}$.
\end{lem}
\begin{proof}
We will show more generally that for a link $L=L_1\cup L_2$, $V(L_1\cup L_2^2)$ is
divisible by $(\xi-1)^{\sharp(L_2)}$, where $L_2^2$ is the two-parallel of $L_2$
with respect to the zero-framing and $\sharp(L_2)$ is the number of components in $L_2$.
\par
We will use L.H.~Kauffman's bracket polynomial \cite{Kauffman:TOPOL87} $\lang\cL\rang$
 for a framed link $\cL$ defined by
\begin{align*}
  &\lang\posi\rang=A\lang\zero\rang+A^{-1}\lang\infi\rang
  \\
  \intertext{and}
  &\lang\cO\rang=1,
\end{align*}
where $\cO$ is the trivial knot with zero-framing.
Since $V(L)=\bigl.(-1)^{\sharp(L)+1}\lang\cL\rang\Bigr\vert_{A=s^{-1/2}}$,
we only have to prove that $\lang\cL_1\cup\cL_2^2\rang$ is divisible by
 $(A^2+A^{-2})^{\sharp(\cL_2)}$ (note that 
$s+s^{-1}=\sqrt{-1}^{\varepsilon}(\xi^{-4\varepsilon u}-\xi^{4\varepsilon u})$).
Here $\cL=\cL_1\cup\cL_2^2$ is the zero-framed link obtained from $L=L_1\cup L_2^2$.
\par
By simple calculations we have
\begin{align*}
    &\lang\posneg\rang-\lang\negpos\rang
  \\[5mm]
  &\quad
  =(A^2-A^{-2})(A^2+A^{-2})
  \left\{
     \lang\nulinf\rang-\lang\infnul\rang
  \right\},
\end{align*}
where dotted lines indicate the connectivity.
By induction on $\sharp(\cL_2)$, we see that the divisibility of $\lang\cL_1\cup\cL_2^2\rang$
does not change by double crossing change between a single strand and a double strand.
This shows that the divisibility of $\lang\cL_1\cup\cL_2^2\rang$ is greater than or
equal to that of $\lang\cL_1\rang\lang\cL_2^2\rang$ which is greater than or equal to
that of $\lang\cL_2^2\rang$, which is greater than or equal to $\sharp(\cL_2)$
by \cite[Lemma 4.1]{Murakami:1998}.
This completes the proof.
\end{proof}
\begin{proof}[Proof of Theorem~\ref{thm:main}]
We only prove the case where $M$ is obtained from an algebraically split, framed
link since the general case follows from T.~Ohtsuki's diagonalizing
lemma \cite[Lemma 2.4]{Ohtsuki:INVEM96}.
\par
Now we assume that $\theta\in H^1(M;\Z/2\Z)$ is given as described before.  
From the lemmas above, we easily see that
$\sqrt{-1}^{\sum_{l=1}^{c}g_l}\alpha(\cL)\Sigma(\cL)$
belongs to $\Z[1/2,\xi]$ and divisible by $(\xi-1)$
(Note that $\sqrt{-1}[2]\in\Z[\xi]$ and divisible by $(\xi-1)$ and that
we need $1/2$ from $\alpha(\cL)$).
\par
Now from \cite{Kirby/Melvin/Zhang:COMMP93},
$\theta\cup\theta\cup\theta\equiv\sum_{l=1}^{c}g_l\mod 2$.
Therefore $\tau_{2p}^{SU(2)}(M,\theta)\in\Z[\xi]$ if
$\theta\cup\theta\cup\theta\equiv 0\mod 2$ and $\sqrt{-1}\Z[\xi]$ otherwise,
completing the proof.
\end{proof}

\bibliographystyle{amsplain}

\providecommand{\bysame}{\leavevmode\hbox to3em{\hrulefill}\thinspace}

\end{document}